\documentclass[11pt]{article}
\usepackage{amsmath, amsfonts, amssymb}
\usepackage{graphicx, amsthm}
\usepackage{cite}
\usepackage{color}

\newtheorem{theorem}{Theorem}[section]
\newtheorem{lemma}[theorem]{Lemma}

\newtheorem{corollary}[theorem]{Corollary}

\theoremstyle{definition}

\newtheorem{example}[theorem]{Example}
\newtheorem{definition}[theorem]{Definition}
\theoremstyle{remark}
\newtheorem{remark}[theorem]{Remark}
\numberwithin{equation}{section}
\def \R{{\bf R}}

\def \gph{\textrm{gph}\,}
\def \dom{\textrm{dom}\,}

\def \lip{\textrm{lip}\,}

\newcommand{\Ls}{\operatornamewithlimits{Limsup}}

\begin{document}

\begin{center}
{\LARGE Sweeping by a tame process}

\bigskip

\textsc{Aris Daniilidis \& Dmitriy Drusvyatskiy}
\end{center}

\bigskip

\noindent \textbf{Abstract.} We show that any semi-algebraic sweeping process admits piecewise absolutely continuous solutions, and any such bounded trajectory must have finite length. Analogous results hold more generally for sweeping processes definable in o-minimal structures. 
This extends previous work on (sub)gradient dynamical systems beyond monotone sweeping sets. 

\bigskip


\vspace{0.6cm}

\textit{Secondary} 34A26, 49J53, 14P10.

\section{Introduction}
\noindent A classical result of \L ojasiewicz shows that any bounded trajectory of the
gradient system 
\begin{equation*}
\dot{x}=-\nabla f(x),
\end{equation*}
with a real-analytic potential function $f$ on $\R^n$, necessarily has finite length
and hence converges to a critical point of $f$. This conclusion can fail if the analyticity is weakened to infinite
differentiability \cite{PalDem82}, for example. The main ingredient of the argument in the analytic case is
the celebrated \L ojasiewicz inequality, which has been subsequently
generalized by Kurdyka \cite{Kurdyka1998} to smooth functions definable in
an o-minimal structure; see \cite{Dries-Miller96} for the relevant
definitions. The authors of \cite{loja,Ioffe-PAMS} further eliminated the smoothness assumptions, showing that any bounded solution $x\colon[0,\eta)\to \mathbf{R}^{n}$ of the subgradient system 
\begin{equation*}
\dot{x}\in-\partial f(x) \qquad \text{a.e. on } [0,\eta)
\end{equation*}
induced by a semi-algebraic function $f$ (or more generally, by one definable
in an o-minimal structure) has finite length and converges to a generalized
critical point of $f$. Here the subdifferential $\partial f$ is meant in any
reasonable sense, such as the limiting subdifferential or the generalized
gradient; see e.g. \cite{RW98}. With the publications of \cite{talweg,Ioffe-2008}, the close relationship between such results and the so-called desingularizing function
(traced back to \cite{Kurdyka1998} for the $C^{1}$ case) became clear. \smallskip

A salient point is that in the subgradient dynamical system, the function $f$
decreases along the trajectory $x$. In particular, after a reparametrization, the trajectory 
$x\colon \lbrack a,b)\rightarrow \mathbf{R}^{n}$ satisfies the inclusion
\begin{equation*}
\dot{x}(r)\in-N_{[f\leq r]}(x(r))\qquad \text{ for a.e. }r\in \lbrack a,b).
\end{equation*}
Here $N_{[f\leq  r]}$ denotes the normal cone to the sublevel set 
$\lbrack f\leq r]:=\{x:f(x)\leq r\}$.
See for example \cite{MP-1991, DLS-2010,DDDL-preprint,italiens, DDL2014} for
this point of view. 
Thus the subgradient system is inherently related to a \textquotedblleft
monotonically evolving sweeping set\textquotedblright \ $r\mapsto \lbrack f\leq r]$. This
observation then naturally motivates investigation of trajectory length
of the more general system 
\begin{equation}
\dot{x}(r)\in-N_{S(r)}(x(r))\qquad \text{ for a.e. }r\in \lbrack a,b),
\label{sp1}
\end{equation}
where $S(r)$ is a subset of $\mathbf{R}^{n}$ varying in time. This dynamical
system is precisely the \emph{sweeping process} of Moreau, well-known in
mathematical mechanics (see e.g. \cite{M77,KM2000}), and which has
recently received much attention \cite{mord_sweep,N_thesis,set_var_sweep,col_sweep,sweep_nonregular}. 
In this paper, much akin to the results of Kurdyka and {\L}ojasiewicz, we prove that
bounded absolutely continuous trajectories of the sweeping process, with a semi-algebraic
set-valued mapping $S\colon \mathbf{R}\rightrightarrows \mathbf{R}^{n}$,
have finite length and therefore must converge to an equilibrium point. We
discuss extensions to the degenerate sweeping process \cite{deg}, and
limitations when the sweeping process is state-dependent in the sense of 
\cite{state}. \smallskip

As a byproduct, we prove a convenient set-valued extension of the projection
formula \cite[Proposition~4]{Lewis-Clarke}, and establish a
``desingularization'' result for semi-algebraic set-valued mappings $S\colon 
\mathbf{R}\rightrightarrows \mathbf{R}^{n}$, generalizing the
Kurdyka-\L ojasiewicz inequality for the sublevel set mapping $r\mapsto[f\leq r]$ of a semi-algebraic function $f$. The desingularization philosophy, combined with \cite{col_exist,ben}, allows us to quickly prove that any locally bounded, semi-algebraic sweeping process always admits piecewise absolutely continuous solutions. The overall trend of the arguments
follows along the lines of \cite{loja,Kurdyka1998}, with some important deviations. Nevertheless, we believe that the striking connection of semi-algebraic and o-minimal geometry to the sweeping process, and in particular to nonmonotone evolution equations,
outlined in this paper, will pave the way for new applications and settings to be explored.

The outline of the manuscript is as follows. In Section~\ref{sec:not}, we
record some notation and preliminary results of variational analysis needed
in the rest of the paper. In Section~\ref{sec:semi}, we discuss basic
elements of semi-algebraic geometry and their interactions with variational constructions. Section~\ref{sec:main} contains our main results on
the lengths of trajectories generated by the sweeping process.
Section~\ref{sec:desing} discusses the role of desingularization, while Section~\ref{sec:exist} applies desingularization ideas to deduce existence of piecewise absolutely continuous solutions of the sweeping process.
\section{Notation}

\label{sec:not} In this section, we summarize some basic tools we will use.
We follow closely the terminology and notation of \cite{RW98}. Throughout,
we consider a Euclidean space which we denote by $\mathbf{R}^{n}$, along
with an inner product $\langle \cdot, \cdot \rangle$ and the induced norm $\| \cdot \|$. The closed unit ball will be denoted by $\mathcal{B}$. For any
set $Q$ in $\mathbf{R}^{n}$, we let $\text{cl}\, Q$ and $\text{int}\, Q$
denote the closure and the interior of $Q$ respectively. The symbol $\text{conv}\, Q$ will stand for the convex hull of $Q$, while $\text{par}\, Q$ will
denote the smallest affine space containing $Q$, translated to the origin,
that is, the linear span of the set $Q-Q$. Given two sets $Q$ and $L$, we
say that the orthogonality relation $Q\perp L$ holds, if any pair of points $x\in Q$ and $y\in L$ are orthogonal. The \emph{distance} of a point $x$ to a
set $Q$ is defined by 
\begin{align*}
\text{dist}\,(x;Q):=\inf_{y\in Q} \|x-y\|
\end{align*}

A set-valued mapping $F$ from $\mathbf{R}^{n}$ to $\mathbf{R}^{m}$, denoted $F\colon \mathbf{R}^{n}\rightrightarrows \mathbf{R}^{m}$, is a mapping from 
$\mathbf{R}^{n}$ to the powerset of $\mathbf{R}^{m}$. The \emph{domain} and \emph{graph} of such a mapping $F$ are defined by 
\begin{align*}
\text{dom}\, F & :=\{x\in \mathbf{R}^{n}: F(x)\neq \emptyset \}, \\
\text{gph}\, F & :=\{(x,y)\in \mathbf{R}^{n}\times \mathbf{R}^{m}: y\in
F(x)\},
\end{align*}
respectively. The inverse of a set-valued mapping $F$ is another set-valued mapping defined by $F^{-1}(y):=\{x:y\in F(x)\}$.
A set-valued mapping $L\colon \mathbf{R}^{n}\rightrightarrows \mathbf{R}^{m}$
is \emph{positively homogeneous} whenever $\text{gph}\, L$ is a cone, or
equivalently whenever we have 
\begin{equation*}
0\in L(0) \qquad \text{and} \qquad L(\lambda x)=\lambda L(x)\quad \text{for
	all } \lambda>0 \text{ and } x\in \mathbf{R}^{n}.
\end{equation*}
In this case, the {\em outer norm} of $L$ is defined by 
\begin{equation*}
|L|^{+} :=\sup_{x\in \mathcal{B}}\, \sup_{y\in L(x)} \|y\|.
\end{equation*}
Due to positive homogeneity of $L$, the outer norm coincides with
\begin{equation*}
\inf \{\kappa>0: \|y\| \leq \kappa\|x\| \textrm{ whenever } y\in L(x)\}.
\end{equation*}
One can now easily deduce that the norm of the inverse admits
the representation 
\begin{equation}  \label{eq:inv}
|L^{-1}|^{+}=\frac{1}{\inf_{\|x\|=1}\text{dist}\,(0;L(x))}.
\end{equation}
Next we pass to certain geometric constructions associated to sets in $\mathbf{R}^{n}$. In what follows, the symbol ``$o(\|x-\bar{x}\|)$ for $x\in Q$'' will stand for any function satisfying $\frac{o(\|x-\bar{x}\|)}{\|x-\bar{x}\|}\to 0$ as $x$ tends to $\bar{x}$ in $Q$.

\begin{definition}[Normal cones]
	Consider a set $Q\subset \mathbf{R}^{n}$ and a point $\bar{x}\in Q$. Then
	the \emph{Fr\'{e}chet normal cone} to $Q$ at $\bar{x}$, denoted $\hat{N}_{Q}(\bar{x})$, consists of all vectors $v\in \mathbf{R}^{n}$ satisfying
	\begin{equation*}
	\langle v, x-\bar{x}\rangle \leq o(\|x-\bar{x}\|) \quad \text{ for } x\in Q.
	\end{equation*}
	The \emph{limiting normal cone} to $Q$ at $\bar{x}$, denoted by $N_{Q}(\bar {x})$, consists of all vectors $v\in \mathbf{R}^{n}$ such that there exist
	sequences $x_{i}$ in $Q$ and $v_{i}\in \hat{N}_{Q}(x_{i})$ satisfying $x_{i}\to \bar{x}$ and $v_{i}\to v$. The \emph{Clarke normal cone} to $Q$ at 
	$\bar {x}$ is simply the set $N^{c}_{Q}(\bar{x}):=\text{cl}\, \text{conv}\,
	N_{Q}(\bar{x})$.
\end{definition}
When $Q$ is a closed convex set, the three normal cones all coincide with the usual convex cone of convex analysis, while for a $C^1$-smooth manifold $Q$ the normal cones coincide with normal spaces in the sense of differential geometry.

Normal cones to graphs of set-valued mappings $F\colon \mathbf{R}^{n}\rightrightarrows \mathbf{R}^{m}$ are naturally associated with
generalized differentiation. Here, we should note that in general the
limiting normal cone $N_{\text{gph}\, F}$ provides much finer information
about the local behavior of $F$ as opposed to the convexified construction $N^{c}_{\text{gph}\, F}$. On the other hand, the results in this paper are
strong enough to pertain to the latter, and hence that's the one we mostly
focus on. Analogous results for limiting constructions are direct
consequences.

\begin{definition}[Coderivatives]{\hfill \\ }
	Consider a mapping $F\colon \mathbf{R}^{n}\rightrightarrows \mathbf{R}^{m}$
	and a pair $(\bar{x},\bar{y})\in \text{gph}\, F$. The \emph{Clarke
		coderivative of }$F$ \emph{at} $\bar{x}$ \emph{for} $\bar{y}$ is the
	set-valued map $D^{*}_{c} F(\bar{x}|\bar{y})\colon \mathbf{R}^{m}
	\rightrightarrows \mathbf{R}^{n}$ defined by 
	\begin{equation*}
	D^{*}_{c} F(\bar x|\bar y)(u):=\{v: (v,-u)\in N^{c}_{\text{gph}\, F}(\bar x,\bar y)\}.
	\end{equation*}
	The \emph{limiting coderivative} $D^{*} F(\bar{x}|\bar{y})$ is defined
	analogously.
\end{definition}

When $F\colon \mathbf{R}^{n}\to \mathbf{R}^{m}$ is $C^{1}$-smooth, then in
terms of $\bar{y}:=F(\bar{x})$, the coderivative mapping $D^{*} F(\bar{x}|\bar{y})$ is single-valued and linear, and coincides with the adjoint of the Jacobian $\nabla F(\bar{x})^{*}$.
Analogously to the smooth case, we use the following notation.

\begin{definition}[Criticality]
	Given a set-valued mapping $F\colon \mathbf{R}^{n}\rightrightarrows \mathbf{R}^{m}$, we say that a pair $(\bar{x},\bar{y})$ in the graph $\text{gph}\, F$
	is a \emph{Clarke critical pair} whenever 
	\begin{equation*}
	0\in D^{*}_{c} F(\bar{x}|\bar{y})(u) \quad \text{ for some } u\neq0.
	\end{equation*}
	A vector $\bar{y}\in \mathbf{R}^{m}$ is a \emph{Clarke critical value} of $F$
	if there exists a point $\bar{x}\in \mathbf{R}^{n}$ so that the pair $(\bar x,\bar{y})$ is Clarke critical. 
\end{definition}

In general, the coderivatives $D^{*}F(\bar{x}|\bar{y})$ and $D^{*}_{c}F(\bar x|\bar{y})$ are positively homogeneous. Hence in particular they admit an
outer norm. Unwrapping the notation for ease of reference, we have 
\begin{equation*}
|D^{*}_{c} F(\bar{x}|\bar{y})|^{+}=\sup_{\|u\| \leq1}\sup \, \{ \|v\|: v\in
D^{*}_{c} F(\bar{x}|\bar{y})(u) \},
\end{equation*}
and 
\begin{equation*}
|D^{*}_{c} F^{-1}(\bar{y}|\bar{x})|^{+} =\frac{1}{\inf_{\|u\|=1} \{
	\|v\|:v\in D^{*}_{c} F(\bar{x}|\bar{y})(u)\}}.
\end{equation*}
In particular $(\bar{x},\bar{y})$ is Clarke critical if and only if $|D^{*}_{c} F^{-1}(\bar{y}|\bar{x})|^{+}=\infty$. \smallskip

\begin{definition}[Asymptotic critical values]
	Given a set-valued mapping $F\colon \mathbf{R}^{n}\rightrightarrows \mathbf{R}^{m}$, we say that a vector $\bar{y}$ is an {\em asymptotic Clarke critical value of} $F$ {\em on a set} $\mathcal{U}\subset\R^m$ if there exists a sequence $(x_i,y_i)\in \gph F$ with $x_i\in \mathcal{U}$, such that $y_i$ converges to $\bar{y}$ and $|D^{*}_{c} F^{-1}(y_i|x_i)|^{+}$ tends to infinity. 
\end{definition}

It is important to note that the outer norm of the limiting coderivative is
very closely related to a pseudo-Lipschitz modulus of the mapping, which will play an important role in Section~\ref{sec:exist}.

\begin{definition}[Aubin Property]
	A set-valued mapping $F\colon \mathbf{R}^{n}\rightrightarrows \mathbf{R}^{m} $ has the \emph{Aubin property at} $\bar{x}$ {\em for} $\bar{y}\in F(\bar{x})$
	if the graph $\text{gph}\, F$ is locally closed around $(\bar {x},\bar{y})$,
	and there are neighborhoods $X$ of $\bar{x}$ and $Y$ of $\bar{y}$, along
	with a constant $\kappa \geq 0$ such that 
	\begin{equation*}
	F(x^{\prime})\cap Y\subset F(x)+\kappa \|x^{\prime}-x\| \mathcal{B}\,,\qquad 
	\text{ for all }\, x,x^{\prime}\in X.
	\end{equation*}
	The infimum of $\kappa$ over all combinations $\kappa$, $X$, and $Y$ so that
	the condition above holds is the \emph{Lipschitz modulus of} $F$ \emph{at} $\bar{x}$ \emph{for} $\bar{y}$, and is denoted by $\text{lip}\, F(\bar{x}|\bar{y})$.
\end{definition}

Provided that the graph of $F$ is closed around $(\bar{x},\bar{y})$,
the following relationships hold: 
\begin{equation*}
\text{lip}\, F(\bar{x}|\bar{y})= \big|D^{*}F(\bar{x}|\bar{y})\big|^{+}\leq 
\big|D^{*}_{c} F(\bar{x}|\bar{y})\big|^{+}.
\end{equation*}
The first equality is a celebrated characterization of the Aubin property
(see \cite{Ioffe-2000}, \cite{boris-book} or \cite[Theorem~9.40]{RW98}, for
example), while the last inequality is immediate from coderivative
definitions.

\smallskip

\section{Semi-algebraic and o-minimal geometry}

\label{sec:semi} A \emph{semi-algebraic} set $Q\subset \mathbf{R}^{n}$ is a
finite union of sets of the form 
\begin{equation*}
\{x\in \mathbf{R}^{n}: f_{1}(x),\ldots,f_{k}(x)=0,g_{1}(x)<0,\ldots
,g_{l}(x)<0\},
\end{equation*}
where $f_{1},\ldots,f_{k}$ and $g_{1},\ldots,g_{l}$ are real polynomials in $n$ variables. It follows immediately that the class of semi-algebraic sets
is closed under the standard Boolean operations (finite unions/intersections
and complementary), while the famous \emph{Tarski--Seidenberg principle} --
also known as \emph{quantifier elimination} -- shows that
semi-algebraicity is preserved under projections. \smallskip

Semi-algebraic subsets of the real line $\mathbf{R}$ are exactly the finite
unions of intervals. This property, known as the \emph{o-minimal}
(order-minimal) property, is the basis for an elegant axiomatization of
various favorable properties of semi-algebraic sets, culminating with a notion of \emph{definable sets}, or more formally, {\em sets definable in an o-minimal structure} \cite{Dries-Miller96}. This
theory allows consideration of much more general sets such as those that are \emph{globally subanalytic}, or sets belonging to the \textrm{log-exp} structure. A slightly more general notion is that of a \emph{tame set} -- a set whose intersection with any ball is definable in an o-minimal structure. Typical examples of tame sets which are not definable, are graphs of \emph{real-analytic} functions. Tame sets are the context of the current paper. We do not however
formally state what definable and tame sets are here since it would take us far off-field. Indeed, the reader can safely replace the word tame (or definable) by semi-algebraic, throughout. We point the interested reader to
the manuscript \cite{Dries-Miller96} or to the short discussion in \cite[p.~771]{Kurdyka1998}.

A key property of definable sets is that they can always be decomposed into a disjoint finite union of smooth (to an arbitrary order) manifolds that fit together in a regular pattern. 
In what follows $p$ will always denote a finite integer no smaller than one.

\begin{definition}[Whitney (a)-regular $C^{p}$-stratification] {\hfill \\ }
	A \emph{Whitney (a)-regular $C^{p}$-stratification} of a set $Q\subset 
	\mathbf{R}^{n}$ is a partition of $Q$ into finitely many $C^{p}$ manifolds
	(called \emph{strata}) satisfying the following compatibility conditions:
	\begin{description}
		\item[\quad Frontier condition:] For any two strata $L$ and $M$, the
		implication 
		\begin{equation*}
		L\cap \text{cl}\,M\neq \emptyset \quad \Longrightarrow \quad L\subset \left( 
		\text{cl}\,M\right) \setminus M\qquad \text{ holds}.
		\end{equation*}
		
		\item[\quad Whitney condition (a):] For any sequence of points $x_{i}$ in a
		stratum $M$ converging to a point $\bar{x}$ in a stratum $L$, if some
		corresponding normal vectors $v_{i}\in N_{M}(x_{i})$ converge to a vector $v$, then the inclusion $v\in N_{L}(\bar{x})$ holds.
	\end{description}
\end{definition}

Definable sets always admit Whitney (a)-regular $C^{p}$-stratifications for
any finite $p$. The importance of such a result for variational analysts can
already be appreciated by observing that the normal cone $N^{c}_{Q}(x)$ must
be contained in the normal space $N_{M}(x)$, where $M$ is a stratum
containing $x$ in any Whitney (a)-regular $C^{1}$-stratifications of $Q$. We
refer the reader to \cite{Lewis-Clarke,Ioffe-PAMS,gen} for applications of
this fact, and of stratifications more broadly, in Variational Analysis. The
forthcoming Theorems~\ref{thm:proj} and \ref{sard} are in the same spirit.
\smallskip

A set-valued mapping $F\colon \mathbf{R}^{n}\rightrightarrows \mathbf{R}^{m}$
(respectively, a function $f\colon \mathbf{R}^{n}\rightarrow \mathbf{R}^{m}$
) is called \emph{definable} if its graph $\text{gph}\,F$ is definable. For
instance, the functions $|x-y|$ and $\sqrt{x^{2}+y^{4}}$ are semi-algebraic,
the function $x\,(\sin x)^{-1}$, for $x\in(0,\pi)$, is globally subanalytic,
while the function 
\begin{equation*}
x\mapsto \exp(\sqrt{|x|})\log(|x|+1)
\end{equation*}
is definable in the log-exp structure. \smallskip

\smallskip

The following is a convenient generalization of the ``projection formula'' \cite[Proposition~4]{Lewis-Clarke} to the coderivative setting. 
Henceforth, we use the symbol $\pi_{x}$ to denote the coordinate projection $(x,y)\mapsto x$, and the symbol $T_{\mathcal{M}}(x)$ to denote the tangent
space to a $C^{p}$-manifold $\mathcal{M}$ at $x$.

\begin{theorem}[Extended projection formula]
	\label{thm:proj} Consider a set-valued mapping $F\colon \mathbf{R}^{n}\rightrightarrows \mathbf{R}^{m}$ and a Whitney (a)-regular $C^{1}$-stratification $\{ \mathcal{M}_{i}\}$ of the graph $\text{gph}\, F\subset 
	\mathbf{R}^{n}\times \mathbf{R}^{m}$. Then for any pair $(\bar{x},\bar{y})\in \text{gph}\, F$ in a stratum $\mathcal{M}_{i}$, the orthogonality
	relation holds: 
	\begin{equation*}
	\text{par}\, D^{*}_{c} F(\bar{x}|\bar{y})(u) \perp \pi_{x}\Big( T_{\mathcal{M}_{i}}(\bar{x},\bar{y})\Big).
	\end{equation*}
\end{theorem}

\noindent \textbf{Proof.} Suppose without loss of generality that $D^*_c F(\bar{x}|\bar{y})(u)$ is nonempty. Then by definition of the coderivative we
have the chain of implications 
\begin{equation*}
v\in D^*_c F(\bar{x}|\bar{y})(u)\quad \Longleftrightarrow \quad (v,-u)\in
N^c_{\text{gph}\, F}(\bar{x},\bar{y})\quad \Longrightarrow \quad (v,-u)\in
N_{\mathcal{M}_i}(\bar{x},\bar{y}).
\end{equation*}
Hence for any tangent vector $(a,b)\in T_{\mathcal{M}_i}(\bar{x},\bar{y})\subset \mathbf{R}^n\times \mathbf{R}^m$, we deduce 
\begin{equation*}
\langle D^*_c F(\bar{x}|\bar{y})(u),a\rangle = \langle u,b\rangle,
\end{equation*}
and consequently $\langle D^*_c F(\bar{x}|\bar{y})(u) -D^*_c F(\bar{x}|\bar{y})(u),a\rangle=0$, as claimed. $\hfill \Box$

\bigskip

The \emph{dimension} of a definable set $Q$ is the maximal dimension of any stratum in any Whitney (a)-regular $C^{1}$-stratification of $Q$. It is well-known that this definition is independent of the choice of the stratification. In particular, a definable subset of $\mathbf{R}^{m}$ has measure zero if and only if it has dimension at most $m-1$. The following is analogous to the main result of \cite{Ioffe-PAMS}.

\begin{theorem}[Sard]
	\label{sard} Consider a definable set-valued
	mapping $F\colon \mathbf{R}^{n}\rightrightarrows \mathbf{R}^{m}$ with a closed graph. Then the set of Clarke critical values of $F$ is a
	definable set of dimension at most $m-1$. Moreover, the set of asymptotic Clarke critical values of $F$ on any bounded definable set $\mathcal{U}$ is also a definable set of dimension at most $m-1$.
\end{theorem}

\noindent \textbf{Proof.} The fact that the set of Clarke critical values is definable follows by standard quantifier elimination. Consider now any
Whitney (a)-regular $C^p$-stratification of $\text{gph}\, F$. Suppose the
equa\-li\-ty $|D^*_c F^{-1}(\bar{y}|\bar{x})|^+=\infty$ holds, or equivalently we have $(0,u)\in N^c_{\gph F}(\bar{x},\bar{y})$ for some vector $u\neq 0$. 
Letting $\mathcal{M}_i$ be the manifold containing $(\bar{x},\bar{y})$, we deduce $(0,u)\in N_{\mathcal{M}_i}(\bar{x},\bar{y})$. Hence $\bar{y}$ is a critical
value in the classical sense of the projection $\pi_y:\mathcal{M}_i\to \mathbf{R}^m$. Applying the standard smooth Sard's theorem with $p$
sufficiently large, we deduce that such critical values $\bar{y}$ have measure zero, as claimed. 

Fix now a bounded definable set $\mathcal{U}$ and suppose that there exists a sequence $(x_i,y_i)\in \gph F$ with $x_i\in \mathcal{U}$, such that $y_i$ converges to $\bar{y}$ and $|D^{*}_{c} F^{-1}(y_i|x_i)|^{+}$ tends to infinity. Hence inclusions
$(w_i,u_i)\in N^c_{\gph F}(x_i,y_i)$ hold for some vectors $w_i$, $u_i$ satisfying $\|u_i\|=1$ and $w_i\to 0$. Since $\mathcal{U}$ is bounded and $\gph F$ is closed, we may suppose $x_i\to\bar{x}$ for some point $\bar{x}$ satisfying $(\bar{x},\bar{y})\in\gph F$.
Letting $\mathcal{M}_i$ be a stratum containing $(\bar{x},\bar{y})$ and passing to the limit we deduce $(0,u)\in N_{\mathcal{M}_i}(\bar{x},\bar{y})$ for some $u\neq 0$. Again applying the standard smooth Sard's theorem with $p$
sufficiently large, we deduce that such values $\bar{y}$ have measure zero, thereby completing the proof. $\hfill \Box$
\smallskip

\begin{remark}[Boundedness of $\mathcal{U}$]
	Boundedness of the set $\mathcal{U}$ is necessary for the theorem above to hold. This can be immediately seen even in the single-valued setting. Indeed, following \cite[Page~776]{Kurdyka1998}, define $F(x,y)=\{r: r\geq \frac{x}{y}\}$ and $\mathcal{U}=\{(x,y): y >0\}$. It follows easily that every $r>0$ is an asymptotic critical value of $F$.
\end{remark}

\section{Finite length of bounded trajectories}
\label{sec:main} Our focus is on the trajectory length of the classical \emph{sweeping process}, introduced by Moreau \cite{M77}. Given a set-valued
mapping $S\colon \mathbf{R}\rightrightarrows \mathbf{R}^{n}$, called the sweeping set, we consider absolutely continuous curves $\gamma \colon(a,b)\rightarrow \mathbf{R}^{n}$ 
satisfying the inclusion
\begin{equation}
\dot{\gamma}(r)\in-N_{S(r)}^{c}\big(\gamma(r)\big)\qquad \text{ for a.e. } r\in(a,b). \label{zero}
\end{equation}
See Figure~1 for an illustration.
Moreau's original construction assumed convexity of the sets $S(r)$, in
which case the normal cone $N_{S(r)}^{c}$ becomes the usual normal cone of
convex analysis. Convexity will not play a role in our work, however. 
We first establish the following bound on the speed of the sweeping process.

\begin{figure}[h!]\label{fig:sweep}
	\centering
	\includegraphics[width=0.6\textwidth]{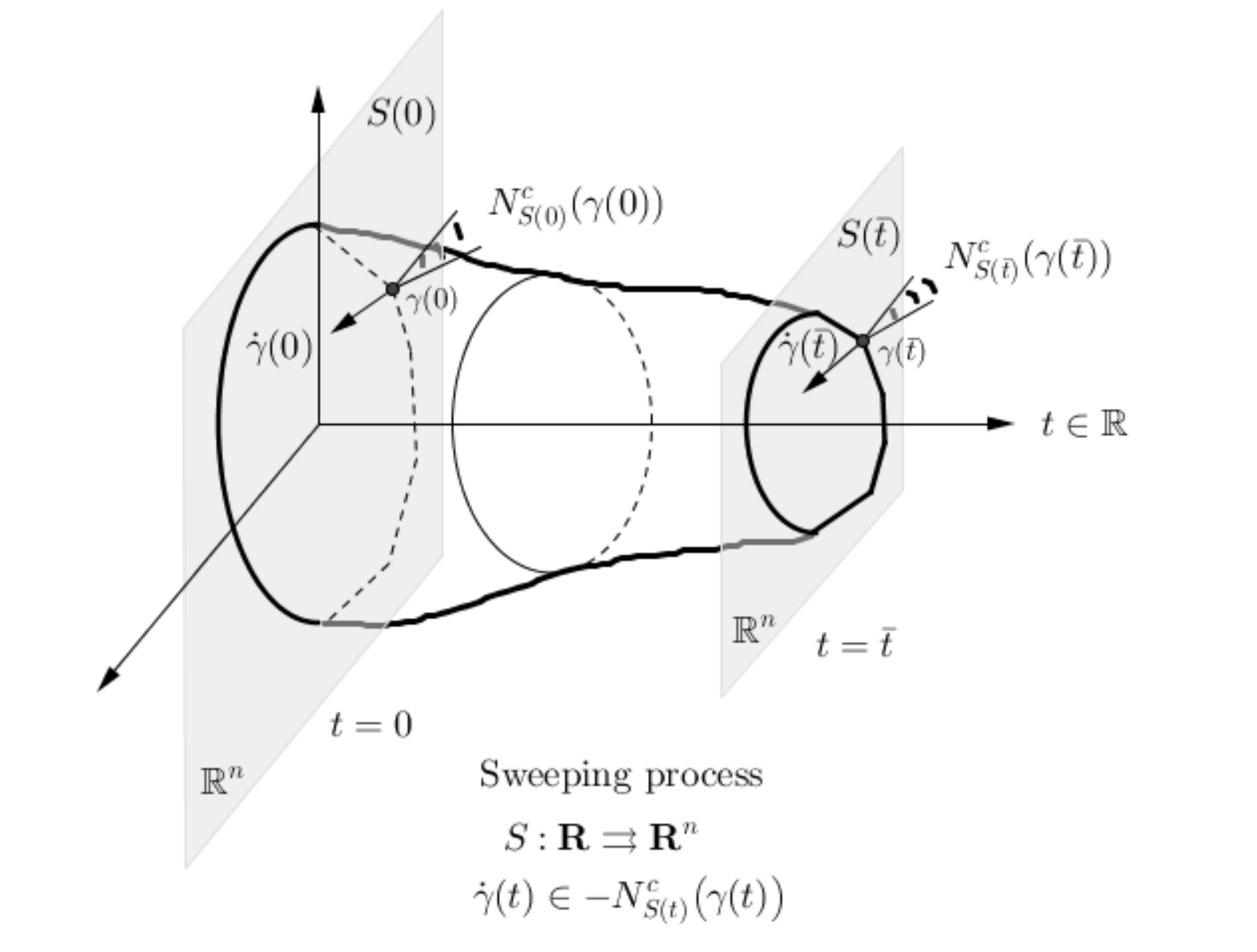}
	\caption{Sweeping process.}
\end{figure}

\begin{theorem}[Speed of the generalized sweeping process]
	\label{Prop_slope} Let $S:\mathbf{R}\rightrightarrows \mathbf{R}^{n}$ be a
	definable mapping with a closed graph and assume that $\gamma \colon
	(a,b)\rightarrow \mathbf{R}$ is \emph{a.e} differentiable and satisfies the
	\textquotedblleft sweeping inclusion\textquotedblright\ \eqref{zero}. Then
	the following estimate on speed holds:
	\begin{equation}
	\Vert \dot{\gamma}(r)\Vert \leq \left\vert D_{c}^{\ast }S\big(r|\gamma (r)
	\big)\right\vert ^{+},\qquad \text{ for a.e. }r\in (a,b).  \label{1}
	\end{equation}
\end{theorem}

\noindent \textbf{Proof.} Notice that \eqref{1} is obvious whenever $\left
\vert D^{\ast}S(r|\gamma(r))\right \vert ^{+}=+\infty,$ thus we may
limit our attention to parameters $r$ with $\left \vert
D^{\ast}S(r|\gamma(r))\right \vert ^{+}<+\infty$ (that is, $S$ has the Aubin
property at $r$ for $\gamma(r)$). Consider the \emph{a.e} differentiable curve
\begin{equation*}
r\mapsto \zeta(r):=(r,\gamma(r))\qquad r\in(a,b),
\end{equation*}
and observe that \eqref{zero} implies $\zeta(r)\in \text{gph}\,S$ for a.e. $r\in(a,b)$.

\medskip

\emph{Claim}. For a.e. $r\in(a,b)$ it holds:
\begin{equation}
\langle N_{\text{gph}\,S}^{c}(\zeta(r)),\dot{\zeta}(r)\rangle =0.
\label{eqn:main_inc}
\end{equation}

\emph{Proof of Claim}\textit{.} Let $\{ \mathcal{M}_{i}\}$ be a Whitney
(a)-regular $C^{1}$-stratification of $\text{gph}\,S$. An easy argument (see
e.g. \cite[Lemma~4.13]{DIL}) shows that for any index $i$ and for a.e. $r\in(a,b)$ the implication 
\begin{equation*}
\zeta(r)\in \mathcal{M}_{i}\quad \Longrightarrow \quad \dot{\zeta}(r)\in T_{\mathcal{M}_{i}}(\zeta(r))\quad \text{ holds}.
\end{equation*}
On the other hand, for such $r$, the Whitney--($a$) condition implies that $N_{\text{gph}\,S}^{c}(\zeta(r))$ is contained in the normal space 
$N_{\mathcal{M}_{i}}(\zeta(r))$. Equation \eqref{eqn:main_inc} follows.$\hfill \lozenge$

\medskip

Fix any $r\in (a,b)$ satisfying \eqref{eqn:main_inc} and assume (with no
loss of generality) that $S$ has the Aubin property at $r$ for $\gamma (r)$.
Setting $H:=\{r\}\times \mathbf{R}^{n}$ we have 
\begin{equation*}
\{r\}\times S(r)=H\cap \text{gph}\,S.
\end{equation*}
Combining this equation with \eqref{zero} we have 
\begin{equation*}
(1,-\dot{\gamma}(r))\in N_{\{r\}\times S(r)}^{c}(\zeta (r))=N_{H\cap \text{gph}\,S}^{c}(\zeta (r)).
\end{equation*}
Since $S$ has the Aubin property at $r$ for $\gamma (r)$, the qualification
condition 
\begin{equation*}
(t,0)\in N_{\text{gph}\,S}(\zeta (r))\quad \Longrightarrow \quad t=0\qquad 
\text{holds}.
\end{equation*}
Hence applying the calculus rule \cite[Theorem 6.42]{RW98}, we
deduce 
\begin{align*}
N_{H\cap \text{gph}\,S}^{c}(\zeta (r))& \subset N_{H}^{c}(\zeta (r))+N_{\text{gph}\,S}^{c}(\zeta (r)) \\
& =\left( \mathbf{R}\times \{0\}\right) +N_{\text{gph}\,S}^{c}(\zeta (r)).
\end{align*}
We conclude that the inclusion 
\begin{equation}
(\lambda ,-\dot{\gamma}(r))\in N_{\text{gph}\,S}^{c}(\zeta (r))  \label{4}
\end{equation}
holds for some $\lambda \in \mathbf{R}$. Appealing then to equation
\eqref{eqn:main_inc}, we obtain the equality 
\begin{equation*}
\langle (\lambda ,-\dot{\gamma}(r)),(1,\dot{\gamma}(r))\rangle =0,\qquad 
\text{ for a.e. }r\in (a,b),
\end{equation*}
and hence 
\begin{equation*}
\lambda =||\dot{\gamma}(r)||^{2}.
\end{equation*}
Plugging this expression back into \eqref{4} we obtain 
\begin{equation*}
\Big(||\dot{\gamma}(r)||^{2},-\dot{\gamma}(r)\Big)\in N_{\text{gph}\,S}^{c}(\zeta (r)).
\end{equation*}
Observe that in the case $\dot{\gamma}(r)=0$, the claimed inequality
\eqref{1} is trivial. Hence we may suppose that this is not the case and
deduce that
\begin{equation*}
\Big(||\dot{\gamma}(r)||,\frac{-\dot{\gamma}(r)}{||\dot{\gamma}(r)||}\Big)
\in N_{\text{gph}\,S}^{c}(\zeta (r)),
\end{equation*}
which readily yields
\begin{equation*}
\left\vert D_{c}^{\ast }S\big(r|\gamma (r)\big)\right\vert ^{+}\geq \Vert 
\dot{\gamma}(r)\Vert ,
\end{equation*}
as claimed.$\hfill \Box $

\bigskip

Following the notation of Theorem~\ref{Prop_slope}, an identical proof shows that
if the stronger inclusion 
\begin{equation*}
\dot{\gamma}(r)\in-N_{S(r)}(\gamma(r))\qquad \text{ holds for a.e. }r\in
\lbrack0,\eta).
\end{equation*}
then the stronger relation 
\begin{equation*}
\Vert \dot{\gamma}(r)\Vert \leq \text{lip}\,S\big(r|\gamma(r)\big).
\end{equation*}
holds for almost every $r$. We note that if $S$ is a Lipschitz continuous mapping, then the authors of \cite{col_exist,ben}  show that the ``catching up algorithm'' produces curves satisfying the above estimate. In contrast, the speed estimate we have just established applies to {\em all} solutions of the sweeping process in the definable setting.

\bigskip

In light of Theorem~\ref{Prop_slope}, to understand the length of the
solutions of the sweeping process it is essential to consider integrability
of the coderivative outer norms $\left \vert D_{c}^{\ast}S\big(r|\gamma (r)
\big)\right \vert ^{+}$. To this end, we introduce the following definition,
much akin to the one considered in \cite[Section~3.3]{talweg} in the context
of subgradient dynamical systems.

\begin{definition}[Talweg function]
	\label{def_talweg}Consider a set-valued mapping $S:(a,b)\rightrightarrows 
	\mathbf{R}^{n}$ and let $\mathcal{U}$ be a subset
	of $\mathbf{R}^{n}$. Then the \emph{talweg function} of $S$ \emph{on} $\mathcal{U}$ is the function 
	$\varphi \colon (a,b)\rightarrow \overline{\mathbf{R}}$ defined by 
	\begin{equation}
	\varphi (r):=\sup_{x\in S(r)\cap \mathcal{U}}\left\{ \left\vert D_{c}^{\ast
	}S\big(r|x\big)\right\vert ^{+}\right\} .  \label{phi}
	\end{equation}
\end{definition}

First, reassuringly the talweg function is rarely infinite.

\begin{lemma}[Finiteness of the talweg function]
	\label{lem:talweg_fin} Consider a definable set-valued mapping $S:(a,b)\rightrightarrows \mathbf{R}^{n}$, with closed values, and let 
	$\mathcal{U}$ be a bounded subset of $\mathbf{R}^{n}$. Then there exists $\epsilon >0$ such that the talweg function of $S$ on $\mathcal{U}$ is not
	equal to $+\infty $ on $(a,a+\epsilon )$.
\end{lemma}

\noindent \textbf{Proof.} 
This follows immediately from Theorem~\ref{sard} applied to $F:=S^{-1}$.
$\hfill \Box $

\bigskip

Next we show that the talwag function of $S$ on a bounded set $\mathcal{U}$ is indeed integrable. The arguments is an application of the curve selection lemma on the
talweg in the spirit of \cite{Kurdyka1998,Ioffe-2008,talweg}.

\begin{theorem}[Integrability of the talweg function]
	\label{Theorem_L1} For any definable, nonempty-valued, closed-valued mapping 
	$S:(a,b)\rightrightarrows \mathbf{R}^{n}$ and the talweg function $\varphi $
	of $S$ on a bounded definable set $\mathcal{U}$, the integral 
	\begin{equation*}
	\int_{a}^{b}\varphi (r)\,dr\quad \text{is finite}.
	\end{equation*}
\end{theorem}

\noindent \textbf{Proof.} For any $r\in (a,b)$, define the set (talweg) 
\begin{equation*}
\mathcal{V}(r):=\left\{ x\in S(r)\cap \mathcal{U}:\left\vert D_{c}^{\ast }S
\big(r|x\big)\right\vert ^{+}\geq \frac{1}{2}\varphi (r)\right\} .
\end{equation*}
Notice that $\varphi $ and $\mathcal{V}$ are definable, and moreover by
Lemma~\ref{lem:talweg_fin} each set $\mathcal{V}(r)$ is nonempty except for
finitely many numbers $r\in (a,b)$. Without loss of generality, assume that
the interval $(a,b)$ does not contain such exceptional points. Using the
curve selection lemma (e.g. \cite[Theorem~1.17]{Dries-Miller96}), we obtain
a definable curve $r\mapsto \theta (r)$ satisfying $\theta (r)\in \mathcal{V}(r)$ for all $r\in (a,b)$. \smallskip

We claim that the equality 
\begin{equation}
t=\Big \langle\dot{\theta}(r),D_{c}^{\ast }S^{-1}(\theta (r)|r)(t)\Big
\rangle\qquad \text{holds for a.e. }r\in (a,b),  \label{eqn:t}
\end{equation}
for which the Clarke coderivative on the right-hand-side is nonempty. To see
this, fix a Whitney (a)-regular $C^{1}$-stratification of $\text{gph}\,S$.
Then for almost every $r$, there exists $\varepsilon >0$ so that the
assignment $\tau \mapsto (\tau ,\theta (\tau ))$ maps the interval $(r-\varepsilon ,r+\varepsilon )$ into a single stratum. For such $r$, fix a
stratum $\mathcal{M}$ containing $(r,\theta (r))$. Then clearly the vector $(1,\dot{\theta}(r))$ is tangent to $\mathcal{M}$ at $(r,\theta (r))$.
Therefore by the Whitney condition (a), the Clarke normal cone $N_{\text{gph}\,S}^{c}(r,\theta (r))$ is contained in the orthogonal complement to 
$(1,\dot{\theta}(r))$. From the definition of the coderivative, we immediately
deduce equality \eqref{eqn:t}, whenever the Clarke coderivative on the
right-hand-side is nonempty. \smallskip

On the other hand, from \eqref{eq:inv} we have unconditionally 
\begin{equation}
\inf_{t=\pm 1}\text{dist}\,(0,D_{c}^{\ast }S^{-1}(\theta (r)|r)(t))=\frac{1}{|D_{c}^{\ast }S(r|\theta (r))|^{+}}.  \label{eq:inverse}
\end{equation}
Hence if neither $-1$ or $+1$ are in the domain of $D_{c}^{\ast
}S^{-1}(\theta (r)|r)$, then we have $0=|D_{c}^{\ast }S(r|\theta
(r))|^{+}\geq \frac{1}{2}\varphi (r)$. On the other hand, for those $r$
satisfying $|D_{c}^{\ast }S(r|\theta (r))|^{+}\neq 0$, equations
\eqref{eqn:t} and \eqref{eq:inverse} immediately imply 
\begin{align*}
1\leq & \inf_{t=\pm 1}\left[ \text{dist}\,(0,D_{c}^{\ast }S^{-1}(\theta
(r)|r)(t))\right] \,\left\Vert \dot{\theta}(r)\right\Vert \\
& =\frac{\left\Vert \dot{\theta}(r)\right\Vert }{|D_{c}^{\ast }S(r|\theta(r))|^{+}}
\leq \frac{2}{\varphi (r)}\left\Vert \dot{\theta}(r)\right\Vert .
\end{align*}
Since the curve $\theta $ is definable and bounded, it has finite length on $(a,b)$, and
consequently the integral $\int_{a}^{b}\varphi (r)\,dr$ is finite, as
claimed. $\hfill \Box $

\bigskip

The main result of the paper is now straightforward.

\begin{corollary}[Bounded length of orbits]
	\label{Cor_length}Let $S:\mathbf{R}\rightrightarrows \mathbf{R}^{n}$ be a
	definable set-valued mapping having a closed graph and let $(a,b)$ be a possibly
	unbounded interval of $\mathbf{R}$. Then any absolutely continuous curve $\gamma \colon (a,b)\rightarrow \mathbf{R}^{n}$ satisfying the sweeping
	inclusion
	\begin{equation*}
	\dot{\gamma}(r)\in -N_{S(r)}^{c}\big(\gamma (r)\big)\qquad \text{ for a.e. }
	r\in (a,b),
	\end{equation*}
	has finite length.
\end{corollary}

\noindent \textbf{Proof.} First, we may assume that the interval $(a,b)$ is
bounded. Indeed, given a trajectory $\gamma :(a,b)\rightarrow \mathbf{R}^{n}$
as above, we consider the semi-algebraic diffeomorphism $\psi :\mathbf{R}
\rightarrow (-1,1)$ by setting $\psi (t)=\frac{t}{\sqrt{1+t^{2}}}$. Then the
curve $\gamma \circ \psi ^{-1}$ is an orbit of the sweeping process $S\circ
\psi ^{-1}$, defined on $(-1,1)$, and it has the same length as $\gamma $.
\smallskip

Without loss of generality, we may also assume that the inclusion $(a,b)\subset \text{dom}\,S$ holds. Since $\gamma $ is bounded, there exists
a bounded set $\mathcal{U}$ containing the image of $\gamma $. Let $\varphi $
be the talweg of $S$ on $\mathcal{U}$. Then by Theorem~\ref{Prop_slope}, for
a.e. $r\in (a,b)$ we have 
\begin{equation}
\Vert \dot{\gamma}(r)\Vert \leq \left\vert D_{c}^{\ast }S\big(r|\gamma (r)
\big)\right\vert ^{+}\leq \varphi (r).  \label{bound_speed}
\end{equation}
Lemma~\ref{Theorem_L1} immediately implies the result. $\hfill \Box $

\bigskip 

\begin{remark}[local monotonicity]
	Theorem \ref{Prop_slope} and Corollary \ref{Cor_length} can be formulated in
	a slightly more general setting, to incorporate certain strongly monotone
	perturbations of the trajectory, as considered in \cite{deg}. To this end,
	recall that a mapping $F\colon \mathbf{R}^{n}\rightarrow \mathbf{R}^{n}$ is 
	\emph{locally $\alpha $-monotone}, whenever we have 
	\begin{equation*}
	\operatornamewithlimits{limsup}_{\Vert x-y\Vert \rightarrow 0}\frac{\langle
		F(x)-F(y),x-y\rangle }{\Vert x-y\Vert ^{2}}\geq \alpha >0.
	\end{equation*}
	Let $\gamma \colon (a,b)\rightarrow \mathbf{R}$ be an absolutely continuous 
	curve and set $\gamma ^{F}:=F\circ \gamma $, where either $F$ or $-F$ is $
	\alpha $-monotone. Assume that $S:\mathbf{R}\rightrightarrows \mathbf{R}^{n}$
	is definable with a closed graph, $\gamma ^{F}$ is absolutely continuous and
	the sweeping inclusion holds:
	\begin{equation*}
	\dot{\gamma}(r)\in -N_{S(r)}^{c}\big(\gamma ^{F}(r)\big)\qquad \text{ for
		a.e. }r\in (a,b).
	\end{equation*}
	Then analogously to Theorem \ref{Prop_slope}, for almost all $r\in (a,b)$ we have
	\begin{equation}
	\Vert \dot{\gamma}(r)\Vert \leq \frac{1}{\alpha }\cdot \left\vert
	D_{c}^{\ast }S\big(r|\gamma ^{F}(r)\big)\right\vert ^{+}\leq \frac{1}{\alpha 
	}\varphi (r),  \label{gen-es}
	\end{equation}%
	where $\varphi $ is the talweg given by (\ref{phi}) for any subset $\mathcal{U}$ of $\mathbf{R}^{n}$ containing the image of 
	$\gamma ^{F}$. 
	The proof of (\ref{gen-es}) follows the lines of the proof of Theorem~\ref{Prop_slope}. Indeed, one considers the curve 
	$r\mapsto \zeta ^{F}(r):=(r,\gamma^{F}(r))$ and eventually deduces
	\begin{equation*}
	\Big(\big\langle\dot{\gamma}(r),\dot{\gamma}^{F}(r)\big\rangle,-\dot{\gamma}
	(r)\Big)\in N_{\text{gph}\,S}^{c}(\zeta ^{F}(r)).
	\end{equation*}
	Considering again the case $\dot{\gamma}(r)\neq 0$ (else the claimed
	inequality is trivial), we get 
	\begin{align*}
	\left\vert D_{c}^{\ast }S\big(r|\gamma ^{F}(r)\big)\right\vert ^{+}& \geq 
	\Big|\Big\langle\frac{\dot{\gamma}(r)}{\Vert \dot{\gamma}(r)\Vert },\dot{\gamma}^{F}(r)\Big\rangle\Big| \\
	& =\lim_{\epsilon \downarrow 0}\Big|\Big\langle\frac{\gamma (r+\epsilon
		)-\gamma (r)}{\Vert \gamma (r+\epsilon )-\gamma (r)\Vert },\frac{(F\circ
		\gamma )(r+\epsilon )-(F\circ \gamma )(r)}{\epsilon }\Big\rangle\Big| \\
	& \geq \alpha \cdot \Vert \dot{\gamma}(r)\Vert ,
	\end{align*}
	and the assertion follows. 
\end{remark}

\bigskip 

It is interesting to note that the analogue of Corollary~\ref{Cor_length}
easily fails when the sweeping set is state-dependent.

\begin{example}[State-dependant process and ODE]
	Consider any autonomous system of ODEs 
	\begin{equation*}
	\dot{x}=F(x),
	\end{equation*}
	where $F$ is a semi-algebraic, Lipschitz continuous mapping. Define the
	semialgebraic set-valued mapping $S(x):=x+F(x)^{\perp}.$ Then every orbit $\gamma$ of the ODE is a solution of the state-dependent sweeping process 
	\begin{equation*}
	\dot{\gamma}(t)\in- N_{S(\gamma(t))}(\gamma(t)).
	\end{equation*}
	Consequently, limit cycles and hence bounded orbits of infinite length can
	easily appear.
\end{example}

\section{Desingularization of the coderivative}\label{sec:desing}

In this section we record a ``desingularization'' result for general definable
set-valued mappings $S:\mathbf{R}\rightrightarrows \mathbf{R}^{n}$ in the
spirit of \cite{Kurdyka1998}. Roughly speaking, any such mapping after a ``localization'' and a reparametrization of its domain can be made to have bounded coderivative norms outside of the critical values of $S^{-1}$.
At the end of the section, we show how our result recovers the desingularization result of Kurdyka \cite{Kurdyka1998}. We use this technique then to investigate solvability of the definable sweeping process in the next section. Here's the main desingularization result.

\begin{theorem}[Desingularization of the Clarke coderivative]
	\label{Thm-main} Consider a definable set-valued mapping $S:\mathbf{R}\rightrightarrows \mathbf{R}^{n}$ having a closed graph, and let $\mathcal{U}
	$ be a bounded subset of $\mathbf{R}^{n}$. Then for any real $a\in \mathbf{R}
	$ there exists a number $b>a$ and a strictly increasing, continuous function 
	$\Psi \colon \lbrack a,b)\rightarrow \mathbf{R}$ that is $C^{1}$-smooth on $(a,b)$, satisfies $\Psi (a)=a$, and such that: 
	\begin{equation*}
	\left\vert D_{c}^{\ast }(S\circ \Psi )(r|x)\right\vert ^{+}\leq 1\qquad 
	\text{ for all }r\in (a,b)\text{ and all }x\in S(\Psi (r))\cap \mathcal{U}.
	\end{equation*}
\end{theorem}

\noindent \textbf{Proof.} If there exists $b>a$ such that the interval $(a,b) $ does not intersect $\text{dom}\,S$, then there is nothing to prove. Consequently, since $S$ is definable, we may suppose that there exists $b>a$
satisfying the inclusion $(a,b)\subset \text{dom}\,S$. Let $\varphi \colon
(a,b)\rightarrow \overline{\mathbf{R}}$ be the talweg of the restriction $S|_{(a,b)}$ on $\mathcal{U}$. Clearly we may assume that $\varphi $ is
continuous on $(a,b)$. If there exists $\varepsilon >0$ such that $\varphi $
equals zero on $(a,a+\epsilon )$, then the theorem is trivial yet again.
Hence we may suppose that $\varphi $ is nonzero on the interval $(a,b)$.
Define now the function 
\begin{equation*}
\Phi (r):=a+\int_{a}^{r}\varphi (\tau )\,d\tau \qquad \text{for }r\in
\lbrack a,b).
\end{equation*}
By Lemma~\ref{Theorem_L1}, the function $\Phi \colon \lbrack a,b)\rightarrow
\lbrack a,\Phi (b))$ above is well defined. Moreover $\Phi $ is clearly
strictly increasing, and $C^{1}$-smooth on $(a,b)$ with a nonvanishing
derivative.

Consider now the inverse $\Psi :=\Phi ^{-1}$. Then $\Psi \colon \lbrack
a,\Phi (b))\rightarrow \lbrack a,b)$ is strictly increasing, conti\-nuous, and 
$C^{1}$-smooth on $(a,\Psi (b))$. Appealing to \cite[Exercise~10.39]{RW98},
for any $\tau \in (a,\Phi (b))$ and any $x\in S(\Psi (\tau ))\cap \mathcal{U}
$ we obtain 
\begin{equation*}
|D_{c}^{\ast }(S\circ \Psi )(\tau |x)|^{+}=\frac{\left\vert D_{c}^{\ast
	}S(\Psi (\tau )|x)\right\vert ^{+}}{\varphi (\Psi (\tau ))}\leq 1,
\end{equation*}
as claimed. $\hfill \Box $

\smallskip
\begin{remark}[Absolute continuity of $\Psi^{-1}$]\label{rem:abs_cont}
	It is immediate from the proof of Theorem~\ref{Thm-main}, that the inverse of the desingularizing function, namely $\Psi^{-1}$, is guaranteed to be absolutely continuous.
\end{remark}

\bigskip 

\subsection{Sweeping by sublevel sets and gradient descent}
We now show how Theorem~\ref{Thm-main} subsumes Kurdyka's seminal desingularization result \cite{Kurdyka1998} for $C^{1}$ definable functions (see also \cite{Lewis-Clarke} for a nonsmooth extension). To this end,
let $f\colon\R^{n}\rightarrow \mathbf{R}$ be a $C^{1}$ definable function
and consider the sweeping process associated to sublevel sets
\begin{equation}
\left\{ 
\begin{array}{l}
S:\mathbf{R}\rightrightarrows \mathbf{R}^{n} \\ 
S(r):=[f\leq r].
\end{array}
\right.   \label{sublevel-map}
\end{equation}
Let $t\mapsto x(t),$ for $t\in \lbrack 0,+\infty )$, be a bounded gradient orbit
for $f$, that is, $\dot{x}=-\nabla f(x)$ with an asymptotic critical value $a:=\lim_{t\rightarrow +\infty }f(x(t)),$ and set $b=f(x(0))$. 
It follows easily that the mapping 
\begin{equation*}
t\mapsto r(t)=b-f(x(t))
\end{equation*}
is a diffeomorphism between $(0,\infty )$ and $(a,b)$. Setting $h=r^{-1}$
and $u=x\circ h$ we obtain a curve $u:(a,b]\rightarrow \mathbf{R}^{n}$ with
the same image as $x$ and satisfying
\begin{equation*}
\dot{u}(r)=-\frac{\nabla f(u(r))}{||\nabla f(u(r))||^{2}},\qquad \text{for }
r\in (a,b].
\end{equation*}
Since equalities $f(u(r))=r$ and $N_{S(r)}(u(r))=\mathbf{R}_{+}\nabla f(u(r))$ hold, we immediately obtain
\begin{equation*}
\dot{u}(r)\in -N_{S(r)}(u(r)).
\end{equation*}
That is, the gradient curve $t\mapsto x(t),$ upon reparametrization, is a
solution of the sweeping process (\ref{sublevel-map}). Moreover, 
an easy computation shows
\begin{equation*}
|D_{c}^{\ast }S(f(x)|x)|^{+}=\frac{1}{||\nabla f(x)||}.
\end{equation*}
Thus the talweg mapping of Definition \ref{def_talweg} reads
\begin{equation*}
\varphi (r)=\left( \inf \,\left\{ \,||\nabla f(x)||:f(x)=r,\,x\in \mathcal{U}
\right\} \right) ^{-1}
\end{equation*}
and Theorem~\ref{Theorem_L1}, Corollary~\ref{Cor_length}, and Theorem~\ref{Thm-main} 
recover the results of Kurdyka in \cite{Kurdyka1998}. 

\section{Existence of solutions}\label{sec:exist}
In this section we will be interested in the existence of trajectories generated by the sweeping process. More specifically, given a set-valued mapping $S\colon [0,\eta)\rightrightarrows\R^n$ and a point $x_0\in S(0)$, we would like to know when there exists a curve $\gamma\colon[0,\eta)\to \R^n$ (appropriately regular) satisfying
\begin{equation}\label{eqn:exist}
\left\{\begin{array}{ll}
-\dot{\gamma}(t)\in N^c_{S(t)}(\gamma(t)) &\quad\textrm{ a.e. on } [0,\eta)\\
\gamma(t)\in S(t) &\quad\textrm{ for all } t\in [0,\eta)\\
\gamma(0)=x_0&
\end{array}\right\}.
\end{equation}
In the case that $S$ is Lipschitz continuous with respect to the Pompeiu--Hausdorff distance, a complete answer was provided in \cite[Theorem 4.2]{col_exist} and \cite[Theorem 3.1]{ben}. Here we mean that a mapping $S\colon [0,\eta)\rightrightarrows\R^n$ is {\em $L$-Lipschitz continuous} if 
$$S(t')\subset S(t)+L|t-t'|\mathcal{B}\qquad \textrm{ for all } t,t'\in [0,\eta).$$
We record below this existence result.

\begin{theorem}[Existence of Lipschitz trajectories]\label{thm:gen_exist}
	Let $S\colon [0,\eta)\rightrightarrows\R^n$ be a $L$-Lipschitz mapping with nonempty, closed values. Then for any $x_0\in S(0)$, there exists a $L$-Lipschitz curve $\gamma\colon[0,\eta)\to \R^n$ satisfying \eqref{eqn:exist}.
\end{theorem}
\smallskip 

\begin{remark}[Extensions to the limiting normal cone]
	In a very recent paper \cite{nad}, it was shown that the analogue of Theorem~\ref{thm:gen_exist} holds for definable $L$-Lipschitz mappings with the limiting normal cone $N_{S(t)}$ replacing the Clarke normal cone $N^c_{S(t)}$. For simplicity, we will state all of our results in the narrower Clarke situation, but an entirely analogous existence theory holds for the limiting case with an identical proof. The only difference is that we must reference the recent manuscript \cite{nad} instead of \cite{col_exist,ben} whenever appropriate.
\end{remark}

There has been a considerable effort recently to weaken the Lipschitz assumption in the theorem above; see for example \cite{ET, CMM} and references therein. 
We will now show that in the definable setting, existence of (at least) piecewise absolutely continuous solutions of \eqref{eqn:exist} can be established even when $S$ is not Lipschitz continuous. This will follow by combining Theorem~\ref{thm:gen_exist} with the desingularization techniques developed in the previous sections.  

We begin with a local existence result. To this end, note that if $S$ is not Lipschitz continuous, then there is an obvious obstruction to having even a continuous local solution of \eqref{eqn:exist} emanating from a point $x_0\in S(0)$. Indeed, when $x_0$ lies outside of the outer limit $$\Ls_{t\searrow 0} S(t):=\left\{\lim_{i\to\infty} x_i: x_i\in S(t_i) \textrm{ with } t_i\searrow 0\right\},$$
clearly no such solution can exist. For example, when $S$ corresponds to a sublevel mapping $S(t)=[f\leq r_0 -t]$ (for a function $f\colon\R^n \rightarrow \R$), such points $x_0\in [f= r_0]\subset S(0)$ are precisely the local minimizers of $f$, and no continuous descent curve can emanate from local minimizers. 
\smallskip
\begin{theorem}[Local existence for the definable sweeping process]\label{thm:exist_loc}
	Consider a definable mapping $S\colon \R_+\rightrightarrows\R^n$ with a closed graph. Then for any 
	$x_0\in \Ls_{t\searrow 0} S(t)$, there exists $\epsilon>0$ and an absolutely continuous curve $\gamma\colon [0,\epsilon)\to \R^n$ satisfying
	\begin{equation*}
	\left\{\begin{array}{ll}
	-\dot{\gamma}(t)\in N^c_{S(t)}(\gamma(t)) &\quad\textrm{ a.e. on } [0,\epsilon)\\
	\gamma(t)\in S(t) &\quad\textrm{ for all } t\in [0,\epsilon)\\
	\gamma(0)=x_0&
	\end{array}\right\}.\end{equation*}
\end{theorem}
\noindent \textbf{Proof.}
Fix a closed ball $\mathcal{U}$ in $\R^n$ containing $x_0$ in its interior, and define the truncation $\widehat{S}(t):=S(t)\cap \mathcal{U}$. Notice, by the assumption $x_0\in \Ls_{t\searrow 0} S(t)$, we have the analogous inclusion $x_0\in \Ls_{t\searrow 0}\widehat{S}(t)$. Appealing to definability, we deduce $(0,\epsilon)\subset \dom \widehat{S}$ for some $\epsilon >0$.  
By Theorem~\ref{Thm-main},  there exists a real number $\eta$ and a strictly increasing, continuous function 
$\Psi \colon \lbrack 0,\eta)\rightarrow \mathbf{R}_+$ that is $C^{1}$-smooth on $(0,\eta)$, satisfies $\Psi (0)=x_0$, and such that: 
\begin{equation*}
\lip(\widehat{S}\circ \Psi )(r|x)\leq 1\qquad 
\text{ for all }r\in (0,\eta)\text{ and all }x\in S(\Psi (r))\cap \mathcal{U}.
\end{equation*}
Shrinking $\eta$, we may assume that the inclusion $(0,\eta)\subset\dom (\widehat{S}\circ \Psi )$ holds.
Appealing to the definition of the Aubin property and the compactness of $\mathcal{U}$, it is easy to see that the mapping $\widehat{S}\circ \Psi $ is locally $1$-Lipschitz continuous around any $r\in (0,\eta)$. Hence  $\widehat{S}\circ \Psi $ is $1$-Lipschitz continuous on the entire interval $(0,\eta)$. Define now the mapping $F\colon [0,\eta)\to\R^n$ given by
\[ F(t):=\begin{cases} 
\widehat{S}\circ \Psi(t) &\quad \textrm{if }t\in (0,\eta) \\
\Ls_{r\searrow 0}\widehat{S}(r) &\quad \textrm{if } t=0  
\end{cases}
\]
Notice that $F(t)$ is 1-Lipschitz continuous, has a closed graph, and satisfies $x_0\in F(0)$. By Theorem~\ref{thm:gen_exist}, then there exist a $1$-Lipschitz curve $x\colon[0,\eta)\to \R^n$ satisfying
\[\left\{\begin{array}{ll}
-\dot{x}(t)\in N^c_{F(t)}(x(t)) &\quad\textrm{ a.e. on } [0,\eta)\\
x(t)\in F(t) &\quad\textrm{ for all } t\in [0,\eta)\\
x(0)=x_0&
\end{array}\right\}.\]
Since the curve $x$ is 1-Lipschitz, shrinking $\eta>0$, we may assume that the image of $x$ is contained in the interior of $\mathcal{U}$. Set $\epsilon:=\lim_{t\uparrow \eta}\Psi(t)$, and define the curve $\gamma\colon [0,\epsilon)\to\R^n$ by setting
$\gamma(r):=x(\Psi^{-1}(r))$. Notice by Remark~\ref{rem:abs_cont}, the inverse $\Psi^{-1}$ is absolutely continuous. Hence $\gamma$ is absolutely continuous as well, being a composition of a Lipschitz function and an absolutely continuous function. Finally observe that $\gamma$ satisfies $\gamma(r)=x(\Psi^{-1}(r))\in F(\Psi^{-1}(r))\subset S(r)$ for all $r\in [0,\epsilon)$ and $-\dot{\gamma}(r)= \frac{-1}{\psi'(\psi^{-1}(r))}\dot{x}(\Psi^{-1}(r))\in N^c_{F(\Psi^{-1}(r))}(x(\Psi^{-1}(r)))=N^c_{S(r)}(\gamma(r))$ 
for a.e. $r\in [0,\epsilon)$. This concludes the proof.
$\hfill \Box$

\smallskip

Next, we will try to maximally extend local solutions of the sweeping process, aiming for a global solution. To this end, we first observe the following.
\smallskip
\begin{corollary}[Convergence to extrema]\label{cor:extrem}
	Consider a definable mapping $S\colon \R_+\rightrightarrows\R^n$ with a closed graph and a point $x_0\in S(t_0)$. Then any absolutely continuous curve $\gamma\colon [0,\epsilon)\to \R^n$, having a maximal domain of definition, such that  
	\[\left\{\begin{array}{ll}
	-\dot{\gamma}(t)\in N^c_{S(t)}(\gamma(t)) &\quad\textrm{ a.e. on } [0,\epsilon)\\
	\gamma(t)\in S(t) &\quad\textrm{ for all } t\in [0,\epsilon)\\
	\gamma(0)=x_0&
	\end{array}\right\},
	\] is either unbounded, or has finite length and converges to some point $x_{\infty}\notin \Ls_{t\searrow \epsilon} S(t)$.
\end{corollary}
\noindent \textbf{Proof.}
This follows immediately from Theorems~\ref{Cor_length} and \ref{thm:exist_loc}. $\hfill \Box$
\bigskip

Let us recall that $S\colon [0,\eta)\rightrightarrows\R^n$ is {\em locally bounded} at $t$ if there exists an open interval $I$ around $t$ such that the image $S(I)$ is a bounded set. We say that $S$ is {\em inner-semicontinuous} at $t$ if for any $x\in S(t)$ and any sequence $t_i\in [0,\eta)$ converging to $t$, there exists a sequence $x_i\in S(t_i)$ converging to $x$. In particular, in the notation of the above theorem, $S$ is not inner-semicontinuous at $\epsilon$, as certified by $x_{\infty}$.  

The following is the main result of the section.

\begin{corollary}[Global existence]{\hfill \\}
	Consider a locally bounded, definable mapping $S\colon [0,\eta)\rightrightarrows\R^n$ with a closed graph and nonempty values. Then for any 
	$x_0\in \Ls_{t\searrow 0} S(t)$ there exists a curve $\gamma\colon [0,\eta)\to\R^n$ satisfying:
	\begin{enumerate}
		\item There is a partition $t_0=0<t_1<\ldots < t_k=\eta$ of the interval $[0,\eta)$ such that $\gamma$ is absolutely continuous on each segment $[t_i,t_{i+1})$; and
		\item The curve $\gamma$ satisfies: 
		\[\left\{\begin{array}{ll}
		-\dot{\gamma}(t)\in N^c_{S(t)}(t) &\quad\textrm{ a.e. on } [0,\eta)\\
		\gamma(t)\in S(t) &\quad\textrm{ for all } t\in [0,\eta)\\
		\gamma(0)=x_0&
		\end{array}\right\}.\] 
	\end{enumerate}
	When $S$ is inner-semicontinuous on the entire interval $[0,\eta)$, then no partition is needed and we can be sure that $\gamma$ is absolutely continuous on the entire interval $[0,\eta)$. 
\end{corollary}
\noindent \textbf{Proof.}
Observe first that since the Aubin property implies inner semicontinuity, by Theorem~\ref{sard} (see also \cite{jeffrey}) the mapping $S$ is inner-semicontinuous at every point $t\in [0,\eta)$ outside of some finite set $\mathcal{N}$. By Theorem~\ref{thm:exist_loc}, there exists $\epsilon >0$ and an absolutely continuous curve $\gamma\colon [0,\epsilon)\to \R^n$ satisfying the conditions \eqref{eqn:exist}. By Zorn's lemma we may maximally extend the domain of $\gamma$ subject to the system \eqref{eqn:exist}. Denote the resulting domain by $[0,a)$. By Corollary~\ref{cor:extrem} and local boundedness of $S$, the curve $\gamma$ converges to some point $x_{\infty}\notin \Ls_{t\searrow a} S(t)$. In particular, $S$ is not inner-semicontinous at $a$ and therefore we deduce that $a\in\mathcal{N}$. We can now repeat the argument with $x_0$ being a point in $\Ls_{t\searrow a} S(t)$. Notice that the latter set is nonempty since $S$ is locally bounded. Concatenating the (finitely many) curves obtained in this way yields the result. $\hfill \Box$

\bigskip

\noindent\textbf{Acknowledgement.} Part of this work has been realized during a research stay of the second author at the University of Chile (November~2013). This author thanks his hosts and the host institution for hospitality. The authors also thank Estibalitz Durand Cartagena and Lionel Thibault for useful discussions.

\bigskip

\noindent Aris Daniilidis

\smallskip

\noindent DIM--CMM, UMI CNRS 2807\newline
Beauchef 851 (Torre Norte, piso~5), Universidad de Chile \smallskip

\noindent E-mail: \texttt{arisd@dim.uchile.cl} \newline
\noindent \texttt{http://www.dim.uchile.cl/{\raise.17ex\hbox{$\scriptstyle\sim$}}arisd}

\smallskip

\noindent Research supported by the grants: \newline
BASAL PFB-03 (Chile), FONDECYT 1130176 (Chile) and MTM2014-59179-C2-1-P
(Spain).\newline

\bigskip

\noindent Dmitriy Drusvyatskiy \newline
University of Washington \newline
Department of Mathematics \newline
C-138 Padelford, Seattle, WA 98195\newline
E-mail: ddrusv@uw.edu \newline
\noindent \texttt{http://www.math.washington.edu/{\raise.17ex\hbox{$\scriptstyle\sim$}}ddrusv/}.
\smallskip

\noindent Research supported by the AFOSR YIP award FA9550-15-1-0237.

\end{document}